\documentclass[reqno]{amsart}
\usepackage{amssymb,amscd}

\begin{document}

\let\kappa=\varkappa
\let\eps=\varepsilon
\let\phi=\varphi

\def\Z{\mathbb Z}
\def\R{\mathbb R}
\def\C{\mathbb C}
\def\Q{\mathbb Q}

\def\OO{\mathcal O}
\def\CP{\C{\mathrm P}}
\def\RP{\R{\mathrm P}}
\def\conj{\overline}
\def\Beta{\mathrm{B}}

\renewcommand{\Im}{{\mathop{\mathrm{Im}}\nolimits}}
\renewcommand{\Re}{{\mathop{\mathrm{Re}}\nolimits}}
\newcommand{\codim}{{\mathop{\mathrm{codim}}\nolimits}}
\newcommand{\id}{{\mathop{\mathrm{id}}\nolimits}}
\newcommand{\Aut}{{\mathop{\mathrm{Aut}}\nolimits}}

\newtheorem{mainthm}{Theorem}
\renewcommand{\themainthm}{\Alph{mainthm}}
\newtheorem{splitthm}{Theorem}[mainthm]
\newtheorem{thm}{Theorem}[section]
\newtheorem{lem}[thm]{Lemma}
\newtheorem{prop}[thm]{Proposition}
\newtheorem{cor}[thm]{Corollary}

\theoremstyle{definition}
\newtheorem{exm}[thm]{Example}
\newtheorem{rem}[thm]{Remark}

\title{Uniformization of strictly pseudoconvex domains}
\author{Stefan Nemirovski}
\address{%
Steklov Mathematical Institute, 119991 Moscow, Russia\hfill\break
\phantom{hh} \& \hfill\break
\phantom{hf} Ruhr-Universit\"at Bochum, D-44780 Bochum, Germany}
\email{stefan@mi.ras.ru}
\author{Rasul Shafikov}
\address{Department of Mathematics, the University of Western Ontario,
London, Canada  N6A 5B7}
\email{shafikov@math.sunysb.edu}
\begin{abstract}
It is shown that two strictly pseudoconvex Stein domains
with real analytic boundaries have biholomorphic universal coverings 
provided that their boundaries are locally biholomorphically equivalent. 
This statement can be regarded as a higher dimensional analogue of the 
Riemann uniformization theorem.
\end{abstract}

\maketitle

\section{Introduction}

Biholomorphic equivalence of domains in complex spaces of dimension
greater than one has proved to be a difficult problem. Already Poincar\'e
observed that the unit ball $B$ in $\C^2$ is not biholomorphically
equivalent to the bidisc $\Delta^2$. In fact, the situation may seem
pretty much hopeless because it has been shown that almost any two randomly
chosen domains in $\C^n$, $n>1$, are inequivalent.

On the positive side, there exist deep results showing that biholomorphic
equivalence of domains with real analytic boundaries is closely
connected to the {\it local\/} biholomorphic equivalence of their boundaries.
For example, Sergey Pinchuk proved that a bounded domain $D\subset \C^n$
with simply connected real analytic strictly pseudoconvex boundary $\partial D$
is equivalent to the unit ball if and only if there is a point in $\partial D$
near which $\partial D$ is equivalent to the unit sphere.

The purpose of this paper is to show how further progress in this direction
can be made by combining the work of Pinchuk~\cite{pi1,pi2,pi3}
on the analytic continuation of biholomorphic maps between real analytic 
hypersurfaces with the classical theory of envelopes of holomorphy
and, more specifically, with the results of Hans Kerner~\cite{Ke2} on the
envelopes of holomorphy of coverings. In this context, it is natural to work
with abstract strictly pseudoconvex domains, not necessarily contained
in~$\C^n$. The general result can be stated as follows:

\begin{mainthm}
\label{mainA}
Let $D$ and $D'$ be Stein strictly pseudoconvex domains with real analytic
boundaries. Then any local equivalence between their boundaries extends
to a biholomorphism of the universal coverings of the {\rm (}open\/{\rm )} 
domains $D$ and~$D'$.
\end{mainthm}

This theorem should in fact be split into two conceptually different cases.
In the generic case, the boundaries are non-spherical, that is, nowhere
locally equivalent to the unit sphere. Then a theorem of Vitushkin 
{\it et al}.~\cite{vek} can be used to obtain a somewhat stronger result:

\begin{splitthm}\label{c2}
If the domains $D$ and $D'$ in Theorem~\ref{mainA} have non-spherical
boundaries, then any local equivalence between $\partial D$ and $\partial D'$
extends to a biholomorphism from the universal covering of $\overline D$ to
the universal covering of $\overline{D'}$.
\end{splitthm}

If the boundaries are spherical, i.\,e., somewhere---and hence, 
by Pinchuk's theorem~\cite{pi1}, everywhere---locally equivalent to the sphere,
Theorem~\ref{mainA} provides the `if' part of the following higher
dimensional analogue of the Riemann uniformization theorem:

\begin{splitthm}\label{c1}
A Stein strictly pseudoconvex domain with real analytic boundary
is uniformized by the unit ball if and only if its boundary is spherical.
\end{splitthm}

In particular, a simply connected Stein domain with compact spherical
boundary is biholomorphic to the ball. This result was proved
by Chern and Ji \cite{cj} for domains $D\subset \C^n$. Recently,
Falbel~\cite{f} observed that it holds also for any simply connected
Stein domain $D$ of complex dimension $n=\dim_\C D>2$.
Indeed, by the Lefschetz theorem for Stein manifolds~\cite{AF},\cite{Bo},
a smoothly bounded Stein domain of complex dimension $n>2$ is
simply connected if and only if its boundary is simply connected,
and therefore this generalization follows from the aforementioned results
of Pinchuk (cp.\ \S\S\,2.3 and~4.2). Note that the question remained
open for general simply connected Stein domains of complex dimension two.

In fact, all the statements above have been proved for domains with simply
connected boundaries by Pinchuk and others (see Section~3). In this
case, the analytic continuation of a local equivalence along the boundary
is single-valued by the monodromy theorem, and the rest is accomplished by
the standard Hartogs theorem. Our main observation is that Kerner's
theorem allows one to extend multiple-valued maps to the domain directly
in a Hartogs-like fashion. Then one can invoke the simple connectivity
of (the universal covering of) the domain and obtain holomorphic 
maps with the desired properties.

Another application of this approach yields a generalization of the result 
of Ivashkovich~\cite{Iv1} on the extension of locally biholomorphic maps 
from real hypersurfaces with non-degenerate indefinite Levi form in~$\CP^n$. 
The r\^ole of the simply connected domain is played here by the complex
projective space itself.

\begin{mainthm}
\label{mainC}
Let $M$ and $M'$ be compact pseudoconcave Levi non-degenerate hypersurfaces 
in $\CP^n$. Suppose that $M$ is real analytic and $M'$ is real algebraic. 
If $M$ and $M'$ are locally equivalent, then the equivalence map is 
the restriction of an automorphism of~$\CP^n$.
\end{mainthm}

Once again, the case of a simply connected $M$ has been treated 
previously in~\cite{hs}. Incidentally, the methods of that paper can be 
used to generalize Theorem~\ref{mainC} to generic pseudoconcave 
CR-submanifolds of higher codimension in~$\CP^n$
with locally injective Segre maps.

\subsection*{Organization of the paper.}
Section 2 concerns the general theory of Riemann domains over
complex manifolds and their envelopes of holomorphy. Kerner's theorem
and its typical applications to analytic continuation are discussed
in~\S 2.3. Then, in \S 2.4 we extend Kerner's theorem to domains over~$\CP^n$
and in \S 2.5 prove a version of the results of Kerner and
Ivashkovich on the extension of locally biholomorphic maps. Section~3
is essentially an overview of known facts about the analytic
continuation of germs of biholomorphic maps between real analytic
hypersurfaces. The last Section~4 contains the proofs of the theorems
stated in the Introduction and a few further corollaries and examples.

\subsection*{Acknowledgments} The authors would like to thank
Sergey Ivashkovich and Seva Shevchishin for helpful discussions.
The first author was partially supported by grants from RFBR
and the program ``Theoretical Mathematics'' of the Russian Academy
of Sciences.


\section{Generalities on analytic continuation}

\subsection{Riemann domains over complex manifolds}
A {\it domain over a complex manifold\/} $X$ is a pair $(D,p)$
consisting of a connected Hausdorff topological space $D$ and
a locally homeomorphic map $p:D\to X$. There exists a unique
complex structure on ~$D$ such that the projection
$p:D\to X$ is a locally biholomorphic map.

A domain $(D,p_D)$ is said to be contained in another
domain $(G,p_G)$ if there is a map $j:D\to G$ such that
$p_G\circ j=p_D$. Notice that the `inclusion' map~$j$
is {\it a~priori\/} only locally biholomorphic and does
not have to be globally injective.

For instance, every ordinary domain ($=$ connected open subset)
$D\subset X$ can be regarded as a domain over $X$ by setting $p_D=\id$.
In this case, the map $j:D\to G$ is a genuine injection because
$p_G\circ j=\id$.

A domain $(D,p_D)$ is called {\it locally Stein\/} if every point $x\in X$
has a neighbourhood $V\ni x$ such that its pre-image $p_D^{-1}(V)\subset D$
is a Stein manifold.

\subsection{Envelopes of holomorphy}
The {\it envelope of holomorphy\/} of a domain $(D,p_D)$ over $X$ is
the maximal domain $(H(D),p_{H(D)})$ over the same manifold~$X$
such that every holomorphic function on $D$ extends to a
holomorphic function on $H(D)$.

More precisely, this means that we are given a locally
biholomorphic map $\alpha=\alpha_D:D\to H(D)$ such that
the following three conditions hold:
\begin{itemize}
\item[(1)] Inclusion: $p_{H(D)}\circ\alpha=p_D$.
\item[(2)]Extension: For every holomorphic function $f\in\mathcal O(D)$
there exists a holomorphic function $F\in\mathcal O(H(D))$ such
that $F\circ\alpha= f$.
\item[(3)] Maximality: If a domain $p_G:G\to X$ and a locally biholomorphic
map $\beta:D\to G$ satisfy conditions (1) and (2) with $G$ in place of $H(D)$,
then there exists a locally biholomorphic map $\gamma:G\to H(D)$ such
that $\gamma\circ\beta=\alpha$ and $p_{H(D)}\circ\gamma=p_G$.
\end{itemize}

The envelope of holomorphy exists and is unique up to a natural isomorphism
by the Thullen theorem. Note once again that the map $\alpha:D\to H(D)$
may not be injective. However, it is injective if $D\subset X$ is a usual
domain in $X$.

Cartan--Thullen and Oka showed that $(H(D),p_{H(D)})$ is a {\it locally\/}
Stein domain over $X$. Oka and Docquier--Grauert proved that every locally
Stein domain over a Stein manifold is Stein. It follows that
{\it the envelope of holomorphy of any domain over a Stein manifold
is Stein}.

Another useful observation is that holomorphic maps to Stein
manifolds extend in the same way as holomorphic functions.

\begin{lem}
\label{SteinVal}
Any holomorphic map $f:D\to Y$ to a Stein manifold extends to a holomorphic
map $F:H(D)\to Y$ {\rm (}in the sense that $F\circ\alpha_D= f$ where
$\alpha_D:D\to H(D)$ is the natural map to the envelope of holomorphy{\rm ).}
\end{lem}

\begin{proof} This is almost obvious. Let $\iota:Y\to\C^N$ be
any proper holomorphic embedding of $Y$ into the complex
vector space of sufficiently high dimension. The components
of the map $\iota\circ f:D\to \C^N$ are holomorphic functions
on $D$ and therefore extend to $H(D)$. The image of the extended
map is contained in $Y=\iota(Y)$ by the uniqueness theorem.\end{proof}

\subsection{Coverings and envelopes over Stein manifolds}
Let $(D,p_D)$ be a domain over a Stein manifold~$X$. Let
$\pi:\widehat D\to D$ be the universal covering of~$D$.
Then the pair $(\widehat D, p_{\widehat D})$ with
$p_{\widehat D}\stackrel{\mathrm{def}}{=}p_D\circ\pi$
is a domain over~$X$. Hence, we can consider its
envelope of holomorphy $(H(\widehat D),p_{H(\widehat D)})$.

\begin{thm}[Kerner~\cite{Ke2}]
\label{Ke1}
The envelope of holomorphy of the universal covering of~$D$
coincides with the universal covering of the envelope of holomorphy of~$D$.
More precisely, there is a commutative diagram of locally biholomorphic
maps\/{\rm :}
$$
\begin{CD}
\widehat D & @>\alpha_{\widehat D}>> & H(\widehat D)\\
@V\pi VV & & @VV H(\pi) V \\
D & @>\alpha_D>> & H(D)
\end{CD}
$$
where the horizontal arrows are the natural maps into
the envelopes of holomorphy and the vertical arrows are the
universal coverings.
\end{thm}

It may be helpful to have in mind the following interpretation
of this theorem in terms of the Weierstra\ss{} theory of analytic
continuation. A holomorphic function on a covering of $D$
corresponds to a germ of a holomorphic function in~$D$
that can be extended analytically along every path in~$D$.
(From this point of view, a covering is a domain over $D$ without
boundary points; see~\cite{GR} for the definition of boundary points.)
Theorem~\ref{Ke1} can be reformulated as follows: {\it If a germ
of a holomorphic function can be extended along every path in $D$,
then it can be extended along every path in the envelope of holomorphy $H(D)$
of $D$ as well.}

\begin{exm}[Monodromy theorem revisited] The classical monodromy
theorem states that if a germ of a holomorphic function can be
extended unboundedly in a simply connected domain, then this
extension is single-valued, i.\,e., defines a holomorphic function.
It follows from Theorem~\ref{Ke1} that the same conclusion holds true if we
assume only that the envelope of holomorphy of the domain is simply
connected.
\end{exm}

\begin{exm}[Pseudoconvex domains]
\label{psconv}
Let $D\subset X$ be a (weakly) pseudoconvex domain with smooth boundary
in a Stein manifold~$X$ of complex dimension $n\ge 2$. Then $D$ is Stein by the
Oka--Docquier--Grauert theorem. Let further $V\subset D$ be the intersection of
a tubular neighbourhood of $\partial D$ in $X$ with the domain~$D$.
It is an open subset of $D$ homotopy equivalent to the boundary~$\partial D$.
Every holomorphic function from $V$ extends to the whole of $D$ by the Hartogs
theorem applied to~$D$. It follows that $V$ is connected and its envelope of
holomorphy is precisely~$D$.

Theorem~\ref{Ke1} shows now that the envelope of holomorphy of
the universal covering of $\pi:\widehat V\to V$ is the universal covering
of the domain~$D$. In other words, if a germ of a holomorphic function
can be extended along every path in $V$, then it extends along every
path in~$D$. In this case, the `enhanced monodromy theorem' asserts
that if $D$ is simply connected, then every holomorphic function
$f\in\OO(\widehat V)$ is the pull-back $F\circ\pi$ of a holomorphic
function~$F\in\OO(D)$.

In view of the latter observation, it is worthwhile to compare the fundamental
groups of $D$ and $V$ or, equivalently, of $D$ and~$\partial D$.
A~standard application of Morse theory to Stein manifolds
(as in~\cite{AF} or~\cite{Bo}) shows that the homomorphism
$\pi_1(V)\to\pi_1(D)$ is an isomorphism if the complex dimension~$n\ge 3$.
Hence, we get no improvement upon the usual monodromy theorem
in this case.

If $n=2$, however, then the homomorphism $\pi_1(V)\to\pi_1(D)$
is only surjective, and the fundamental group of the boundary can be much
larger than that of the domain. For instance, tom Dieck and Petrie~\cite{tDP}
gave explicit examples of affine algebraic surfaces $\Sigma\subset \C^3$
such that the intersection of $\Sigma$ with a sufficiently large ball
is a {\it contractible\/} strictly pseudoconvex domain whose boundary
has {\it infinite\/} fundamental group.
\end{exm}

\subsection{Envelopes over the complex projective space}
Analytic continuation over $\CP^n$ can be quite well understood
with the help of the following theorem of Fujita~\cite{Fu}
and Takeuchi~\cite{Ta}: {\em A locally Stein domain over $\CP^n$
is either Stein or coincides with $\CP^n$ itself}.
An elegant and illuminating proof of this result
was given by Ueda~\cite{Ue}.

In particular, for any domain $(D,p)$ over $\CP^n$, its envelope
of holomorphy $H(D)$ is either a Stein domain over $\CP^n$ or
the tautological domain $(\CP^n,\id)$. Clearly, the latter option
can be characterized by the property that every holomorphic
function on $D$ is constant.

\begin{thm}
Theorem~$\ref{Ke1}$ holds true for
every domain $(D,p)$ over the complex projective space~$\CP^n$.
\end{thm}

\begin{proof}
If the envelope $H(D)$ is Stein, then we can equip $D$ with
the locally biholomorphic map $\alpha_D:D\to H(D)$
and consider it as a domain over $H(D)$. Consequently,
the result follows directly from Kerner's theorem
in this case.

If $H(D)=\CP^n$, then the analogue of Theorem~\ref{Ke1}
is given by the following proposition.\end{proof}

\begin{prop}
\label{cp1}
Let $(D,p)$ be a domain over $\CP^n$ such that every holomorphic
function on $D$ is constant. Then the same holds true for every
covering $\pi:\widehat D\to D$.
\end{prop}

\begin{proof}
It suffices to consider the case of the universal covering
$\pi:\widehat D\to D$ because every holomorphic function on
any covering of $D$ can be pulled back there.

Let us assume that $H(\widehat D)$ is Stein and seek
a contradiction by imitating Kerner's arguments from~\cite[pp.~127--129]{Ke2}.
Let $\Delta$ be the group of deck transformations of $\widehat D$.
The action of this group on $\widehat D$ extends to a properly
discontinuous free action of its factor group $\widetilde\Delta$ on
$H(\widehat D)$. The quotient $H(\widehat D)/\widetilde\Delta$
is a domain over $\CP^n$ containing $D=\widehat D/\Delta$
and covered by~$H(\widehat D)$. Kerner's Hilfssatz~1 shows that
this domain is {\it locally\/} Stein. However, it cannot be Stein
because it would then follow that there exist non-constant holomorphic
functions on~$D$. Hence, $H(\widehat D)/\widetilde\Delta=\CP^n$ by Fujita's
theorem. Since $\CP^n$ is simply connected, it follows that
$H(\widehat D)=\CP^n$, a contradiction.\end{proof}

In order to consider the extension of meromorphic functions as
well, let us introduce the {\it envelope of meromorphy\/} $(M(D),p_{M(D)})$
of a domain $(D,p)$ over a complex manifold $X$. The definition
is completely analogous to the holomorphic case. Namely,
$(M(D),p_{M(D)})$ is the maximal domain over~$X$ containing $(D,p_D)$
and such that every {\it meromorphic\/} function from $D$ extends
meromorphically to~$M(D)$. The existence and uniqueness of the envelope
of meromorphy follow from a Thullen-type theorem. Levi's theorem on the
extension of meromorphic functions and the ubiquitous theorem of Oka
imply that the envelope of meromorphy is a locally Stein domain over $X$.

One can show that the envelope of meromorphy of any domain over
the complex projective space coincides with its envelope
of holomorphy but we shall only need the following partial result.

\begin{prop}
\label{cp2}
Let $(D,p)$ be a domain over $\CP^n$ such that every holomorphic
function on $D$ is constant. Then every meromorphic function
on $D$ has the form $f\circ p$ for a rational function $f$ on $\CP^n$.
\end{prop}

\begin{proof}
We have to show that the envelope of meromorphy of $(D,p)$
coincides with $\CP^n$. Then every meromorphic function on $D$
can be obtained as the pull-back of a meromorphic function on $\CP^n$
which must be rational by Serre's GAGA principle.

Suppose that the envelope of meromorphy is not $\CP^n$. Since
it is locally Stein, it must be a Stein domain by Fujita's
theorem. But this implies that $(D,p)$ is contained in a
Stein domain and hence admits non-constant holomorphic
functions, a contradiction.
\end{proof}

\begin{exm}[Pseudoconcave hypersurfaces in $\CP^n$]
\label{Concave}
Let $M\subset\CP^n$ be a smooth compact real hypersurface whose
Levi form is non-degenerate and indefinite (i.\,e.,
has positive and negative eigenvalues at each point).
Consider a connected neighbourhood $U\supset M$ and
let $H(U)$ and $M(U)$ be its envelopes of holomorphy
and meromorphy, respectively. Note that neither envelope
is Stein because a hypersurface such as $M$ cannot lie
in a Stein manifold. Hence, $H(U)=M(U)=\CP^n$ by the
Fujita theorem. In other words, every holomorphic
function on $U$ is constant and every meromorphic
function on $U$ is rational.

Furthermore, let $\pi:\widehat U\to U$ be a covering of the
neighbourhood~$U$. Then every holomorphic function on $\widehat U$
is constant by Proposition~\ref{cp1} and therefore every meromorphic
function on $\widehat U$ is the pull-back of a rational function
by Proposition~\ref{cp2}. Thus $H(\widehat U)=M(\widehat U)=\CP^n$ as well.
\end{exm}

\subsection{Extension of locally biholomorphic maps}
The theory of envelopes of holomorphy and meromorphy
recalled in this section allows us to give a unified
treatment (and a slight generalization) of the results
on the holomorphic extension of locally biholomorphic
maps established by Kerner~\cite{Ke1} and Ivashkovich~\cite{Iv1}.

\begin{thm}
\label{LocBiExt}
Let $X$ and $Y$ be complex manifolds each of which
is either Stein or biholomorphic to the complex
projective space~$\CP^n$. Let $(D,p)$ be a domain over $X$
and $(H(D),p_{H(D)})$ be its envelope of holomorphy. Then
every locally biholomorphic map $f:D\to Y$ extends to
a locally biholomorphic map $F:H(D)\to Y$.
\end{thm}

\begin{proof}
Let us begin with a simple general observation. Suppose that
$A\subset H(D)$ is a non-empty complex hypersurface ($=$ a
complex analytic subset of pure codimension one). Then $A$
must intersect the image $\alpha_D(D)$ of the domain $D$ in
its envelope of holomorphy. Indeed, the complement $H(D)\setminus A$
is a proper locally Stein open subset of $H(D)$ and hence must
be Stein by the theorems of Oka--Docquier--Grauert (if $H(D)$
is Stein) and Fujita (if $H(D)=\CP^n$). It follows that
there exists a holomorphic function, say, $g\in\OO(H(D)\setminus A)$
which cannot be extended to $H(D)$. If $\alpha_D(D)\cap A$
were empty, then the function $g\circ\alpha_D\in\OO(D)$ would
not extend to $H(D)$, a contradiction.

Assume now that the target manifold $Y$ is Stein. Then the map
$f:D\to Y$ extends to a holomorphic map $F:H(D)\to Y$ by Lemma~\ref{SteinVal}.
Consider the ramification locus of~$F$, i.\,e., the complex
hypersurface in $H(D)$ consisting of the points at which
$\mathop{\mathrm{rank}}\nolimits_\C(F)<n=\dim_\C X=\dim_\C Y$.
This hypersurface cannot intersect $\alpha_D(D)$, and so it is
empty by the observation above. Thus, $F$ is a locally biholomorphic map.

The case $Y=\CP^n$ splits into two subcases. Suppose first that $D$ admits
a non-constant holomorphic function. Since the map $f:D\to\CP^n$ is locally
biholomorphic, we can regard the pair $(D,f)$ as a domain {\it over\/}~$Y$.
Let $(\widetilde D,\widetilde f)$ be the envelope of holomorphy
of this domain over~$Y$. (We use different notation for envelopes
over $Y$ to avoid confusion.) Let $\beta:D\to\widetilde D$ be
the natural map into the envelope. Recall that $\widetilde f\circ\beta=f$.
Since $\OO(D)\ne \C$, the envelope $\widetilde D$ is a Stein
manifold by Fujita's theorem. As we have already seen, the map
$\beta:D\to\widetilde D$ extends to a locally biholomorphic map
$\Beta:H(D)\to\widetilde D$ such that $\Beta\circ\alpha_D=\beta$.
The composition $F\stackrel{\mathrm{def}}{=}\widetilde f\circ\Beta$ is
the desired extension of~$f$ to $H(D)$. Indeed,
$$
F\circ\alpha_D=\widetilde f\circ\Beta\circ\alpha_D=\widetilde f\circ\beta=f.
$$

Finally,  consider the case $\OO(D)=\C$. This is only possible if
$X=\CP^n$ as well. Thus we can apply Proposition~\ref{cp2} to the
components of the map $f$ with respect to an affine coordinate
system on $Y=\CP^n$. It follows that there is a rational map
$F:\CP^n\to\CP^n$ such that $F\circ\alpha_D=f$.
In the same way as before, we see that the ramification locus of $F$
(defined to be the Zariski closure of the set of points at which $F$
is holomorphic and not of maximal rank) must be empty. Therefore,
the map $F$ is locally biholomorphic on the complement of its
indeterminacy locus~$I$. However, $I$ has complex codimension at
least~$2$, so we can apply the preceding case of the theorem
locally, in a Stein neighbourhood of each point of $I$, and conclude
that $F$ is in fact locally biholomorphic everywhere.
\end{proof}

Notice that a locally biholomorphic map from $\CP^n$ to itself is
just a linear automorphism. Hence, we obtain the following
(cp.~\cite{Iv1} and~\cite{hs}):

\begin{cor}
\label{CPRig}
Let $(D,p)$ be a domain over $\CP^n$ such that every holomorphic
function on $D$ is constant. Then every locally biholomorphic
map $f:D\to\CP^n$ has the form $L\circ p$ for an automorphism
$L\in\mathop{\mathrm{Aut}}(\CP^n)$.
\end{cor}

The proof of Theorem~\ref{LocBiExt} rightly suggests that a holomorphic
but not locally biholomorphic map $f:D\to\CP^n$ may not admit a holomorphic
extension to~$H(D)$. For instance, the quadratic transformation $Q:\CP^2\to\CP^2$
is a birational map with three indeterminacy points. Clearly, if $D$ is a
punctured neighbourhood of one of these points, then $Q$ is holomorphic
on $D$ but cannot be holomorphically extended to~$H(D)$.
This phenomenon was discovered by Ivashkovich in~\cite{Iv1}
and became the starting point for his deep results on meromorphic
continuation (see, for instance,~\cite{Iv2}).

\begin{rem} Ueda~\cite{Ue} generalized Fujita's theorem to Grassmann
manifolds $\mathop{\mathrm{Gr}}\nolimits_\C(m,n)$. Hence, the same
proof shows that Theorem~\ref{LocBiExt} and Corollary~\ref{CPRig}
remain valid with $\CP^n$ replaced by any complex Grassmannian.
\end{rem}

\section{Analytic continuation along real hypersurfaces.}

\subsection{Extension of local equivalences between real hypersurfaces}
As it was mentioned in the Introduction, biholomorphic equivalence of
domains in several complex variables is closely related to the equivalence
of their boundaries. In what follows we assume that $n\ge 2$.
Given two real hypersurfaces $M$ and $M'$ in $\C^n$ we say
that they are locally biholomorphically equivalent (or simply equivalent)
at points $p\in M$ and $p'\in M'$ if there exist connected open
neighbourhoods $U\ni p$ and $U'\ni p'$ in $\C^n$ and a biholomorphic
map $f:U\to U'$ with $f(M\cap U)=M'\cap U'$.

It turns out that for real analytic hypersurfaces local equivalence
at a point may imply equivalence everywhere. In this formulation, the first
result was obtained by Pinchuk:

\begin{thm}[\cite{pi1}] 
\label{PiSph}
Let $M\subset \C^n$ be a connected strictly pseudoconvex real analytic
hypersurface which is equivalent to the unit sphere
$S^{2n-1}\subset \C^n$ at some point $p\in M$. Then the germ of the equivalence
map $_p{\bf f}$ extends along any path on $M$ as a locally biholomorphic map
sending $M$ to~$S^{2n-1}$.
\end{thm}

This allows us to define {\it spherical\/} hypersurfaces as being connected
strictly pseudoconvex real analytic and equivalent to $S^{2n-1}$ at some
(and therefore at every) point. It is an immediate corollary of Pinchuk's
result that a compact simply connected spherical hypersurface $M$ in an
arbitrary complex manifold $X$ is biholomorphic to the standard sphere.
Hence, Theorem~\ref{PiSph} can be formulated with the sphere replaced by
any compact simply connected spherical hypersurface $M'$. This result will
be false, however, if $M'$ is only assumed to be spherical and compact
(see~\cite{bs}).

Pinchuk proved in~\cite{pi2} that the compactness assumption suffices
in the {\it non-spherical\/} case, i.\,e., when $S^{2n-1}$ is replaced by
a compact strictly pseudoconvex real analytic hypersurface $M'\subset \C^n$
which is {\it not\/} equivalent to the sphere at any point. This was
generalized to hypersurfaces in arbitrary complex manifolds in~\cite{vek}.
Combining these results gives

\begin{thm}[\cite{pi3},\cite{vek},\cite{Vi}]
\label{NonSph}
Let $X$ and $X'$ be complex manifolds, $\dim_\C X=\dim_\C X'=n$.
Let $M\subset X$ and $M'\subset X'$ be real analytic strictly pseudoconvex
non-spherical hypersurfaces, $M$ connected and $M'$ compact. If $M$ and $M'$
are locally equivalent at points $p\in M$ and $p'\in M'$, then the germ 
$_p{\bf f}$ of the equivalence map extends along any path on $M$ as 
a locally biholomorphic map sending $M$ to~$M'$.
\end{thm}

Much less is known about maps between hypersurfaces which are not 
strictly pseudoconvex. Beloshapka and Ezhov have given examples 
suggesting that an extension of Theorem~\ref{NonSph} to this 
situation may be problematic. On the other hand, the reflection principle 
underlying the proof of Theorem~\ref{PiSph} can still be used in its 
geometric form provided that the target hypersurface is real algebraic 
(cf.~\cite{sh1},\cite{sh2}). For hypersurfaces with non-degenerate 
indefinite Levi form this approach yields the following result: 

\begin{thm}[\cite{hs}]
\label{Shaf}
Let $M$ be a connected real analytic Levi non-degenerate hypersurface 
in a complex manifold $X$, and let $M'$ be a compact real algebraic 
Levi non-degenerate hypersurface in~$\CP^n$. Suppose that $M$ and $M'$ 
are locally equivalent at points $p\in M$ and $p'\in M'$. Then the germ 
$_p{\bf f}$ of the equivalence map extends along any path on $M$ as 
a locally biholomorphic map sending $M$ to~$M'$.
\end{thm}

This theorem was applied in~\cite{hs} to prove Theorem~\ref{mainC} 
for a simply connected hypersurface~$M\subset\CP^n$. Let us also note
that if both $M$ and~$M'$ are real algebraic, Theorem~\ref{Shaf} follows 
essentially from the well-known theorem of Webster~\cite{w}, which states 
that if $M$ and $M'$ are algebraic hypersurfaces in $\C^n$, locally equivalent 
at $p\in M$ and $p'\in M'$, and $p$ and $p'$ are Levi non-degenerate points, 
then the equivalence map is algebraic, i.\,e., its graph is contained in 
an algebraic subvariety of $\C^n \times \C^n$.

\subsection{Analytic continuation to open neighbourhoods}
If $M$ is not simply connected, then the analytic continuation of an
equivalence germ $_p{\bf f}:M\to M'$ along homotopically nonequivalent
paths starting at $p\in M$ with the same end point $q\in M$ may produce
different extensions, which {\it a priori\/} may have different radii
of convergence at~$q$. Nonetheless, an extension into a fixed open set
can be sometimes obtained from general theory.

\begin{lem}
\label{LocExtSt}
Let $M\subset X$ be a compact strictly pseudoconvex hypersurface and
$D'$ a Stein strictly pseudoconvex domain. Suppose that a germ
$_{p}{\bf f}:M\to \partial D'$ of a locally biholomorphic map
extends along any path in~$M$ as a locally biholomorphic map sending
$M$ to $\partial D'$. Then there exists a neighbourhood $U\subset X$
of $M$ such that $_{p}{\bf f}$ extends as a locally biholomorphic map
with values in $D'$ along any path in $U^-$, the strictly pseudoconvex 
one-sided neighbourhood of~$M$. 
\end{lem}

\begin{proof}
It is enough to prove that for any point $q\in M$ there exists a neighbourhood
$U_q\subset X$ such that any holomorphic map $g$ obtained by analytic continuation
of $_{p}{\bf f}$ extends to the strictly pseudoconvex side of $U_q$. 
Let $V$ be a neighbourhood of $q$ such that $M\cap V$ is simply connected. 
Then $g$ extends to a locally biholomorphic map in some neighbourhood 
of $M\cap V$. This extension takes values in $D'$ on the pseudoconvex 
side of~$M$. Since $M$ is strictly pseudoconvex, there exists a neighbourhood 
$U_q\subset X$ of the point~$q$ such that every function holomorphic 
in a one-sided neighbourhood of $M\cap V$ extends holomorphically
to $U_q^-$, the pseudoconvex side of $U_q$. The same holds true for locally 
biholomorphic maps into the Stein manifold~$D'$ by Theorem~\ref{LocBiExt}. 
By construction, $U_q$ is independent of $g$.
\end{proof}

\begin{lem}[\cite{Iv1}]
\label{LocExtCP}
Let $M\subset X$ be a compact Levi non-degenerate pseudoconcave hypersurface.
Suppose that a germ $_{p}{\bf f}:M\to\CP^n$ of a locally biholomorphic map
to the complex projective space extends as a locally biholomorphic map
along any path on $M$. Then there exists a neighbourhood $U\subset X$ of $M$
such that $_{p}{\bf f}$ extends as a locally biholomorphic map along
any path in $U$.
\end{lem}

\begin{proof}
The argument is very similar to the previous one. For a point $q\in M$,
let $V\ni q$ be a coordinate neighbourhood such that $M\cap V$
is simply connected. By the Hans Lewy theorem, there exists a neighbourhood
$U_q\subset V$ of the point~$q$ such that every function holomorphic 
in a neighbourhood of $M\cap V$ extends holomorphically to $U_q$. 
By Theorem~\ref{LocBiExt} the same extension property holds for 
locally biholomorphic maps to $\CP^n$. Hence, any holomorphic map obtained 
by analytic continuation of $_{p}{\bf f}$ extends to a locally biholomorphic 
map $U_q\to\CP^n$.
\end{proof}

\section{Global extension of local maps}

\subsection{General set-up and the proof of Theorem~\ref{mainA} in the non-spherical case}
Let $D$ and $D'$ be two Stein strictly pseudoconvex domains with real analytic
boundaries. Suppose that $\partial D$ and $\partial D'$ are locally equivalent
at points $p\in \partial D$ and $p'\in \partial D'$.
Let $_p{\bf f}:\partial D\to\partial D'$ be the germ of a locally
biholomorphic map realizing this equivalence.

\begin{prop}
\label{global}
If the germ $_p{\bf f}$ can be extended as a locally biholomorphic map
along any path in~$\partial D$, then it can be extended as a locally
biholomorphic map with values in $\conj{D'}$ along any path in~$\conj D$.
\end{prop}

\begin{proof}
Since the boundaries are real analytic, it follows that any map obtained by
the extension of $_p{\bf f}$ along a path in the boundary will send $\partial D$
to~$\partial D'$.

By Lemma~\ref{LocExtSt}, there exists a neighbourhood
$U\supset\partial D$ such that the germ $_p{\bf f}$ can be extended
along any path in $V=U\cap D$ as a locally biholomorphic map with values in~$D'$.
This extension defines a locally biholomorphic map $f:\widehat V\to D'$ from the
universal covering $\widehat V\to V$.

The envelope of holomorphy of $V$ is precisely $D$ by the
Hartogs theorem (cf.\ Example~\ref{psconv}). Hence, the envelope of holomorphy
of the universal covering $\widehat V\to V$ is the universal covering
$\pi:\widehat D\to D$ by Kerner's theorem.
Theorem~\ref{LocBiExt} shows now that the map $f:\widehat V\to D'$ extends
to a locally biholomorphic map $F:\widehat D\to D'$.

Let finally $\conj\pi:\conj Y\to\conj D$ be the universal covering of
the closure of~$D$. Then $\conj Y$ is a complex manifold with (not necessarily
compact) boundary $\partial Y=\conj\pi^{-1}(\partial D)$
and interiour $Y=\conj Y-\partial Y=\widehat D$.
By construction, the map $F:Y\to D'$ coincides with a lift of an extension of $_p{\bf f}$
near every boundary point $q\in\partial Y$. Hence, it extends to a
locally biholomorphic map $\conj F:\conj Y\to\conj{D'}$ of complex manifolds
with boundary. In other words, the germ $_p{\bf f}$ extends as a locally
biholomorphic map with values in $\conj{D'}$ along any path in~$\conj D$.
\end{proof}

We are now in position to prove the stronger form of Theorem~\ref{mainA}
for non-spherical domains (Theorem~\ref{c2}). Indeed, let
$_p{\bf f}:\partial D\to\partial D'$ be a local equivalence germ.
By Theorem~\ref{NonSph}, this germ extends along any path in $\partial D$.
Therefore, it extends as a locally biholomorphic map along any path
in $\conj D$ by Proposition~\ref{global}. The same conclusions hold true
for the inverse map $_{p'}{\bf f}^{-1}:\partial D'\to\partial D$ as well.
Hence, it follows from the monodromy theorem that the extension of $_p{\bf f}$
defines a biholomorphism of the universal covering of the closure of $D$
onto the universal covering of the closure of $D'$.

\begin{cor}
Suppose that $D$ and $D'$ are Stein strictly pseudoconvex domains
with locally equivalent non-spherical boundaries. If $D$ has finite
fundamental group, then so does~$D'$.
\end{cor}

\begin{proof}
The fundamental group of a compact manifold with boundary (e.\,g., $\conj D$) 
is finite if and only if the universal covering of this manifold is compact.
\end{proof}

\subsection{Uniformization of Stein domains with spherical boundary}
Any local equivalence between spherical hypersurfaces factors through the
sphere. Therefore, we only need to prove Theorem~\ref{mainA} in the case
when $D'$ is the unit ball $B\subset\C^n$. Let $S=\partial B$ denote the 
unit sphere. By Theorem~\ref{PiSph} and Proposition~\ref{global},
any local equivalence germ $_p{\bf f}:\partial D\to S$ extends to
a locally biholomorphic map $\conj F:\conj Y\to \conj B$ of complex 
manifolds with boundary from the universal covering of $\conj D$ 
to the closed unit ball. 

The map $\conj F$ may be viewed as the extension of the `developing map'
of the boundary, introduced by Burns--Shnider~\cite{bs}, to the universal 
covering of the domain. In particular, it inherits the following important 
equivariance property:

\begin{lem}
\label{equiv}
Let $\Gamma=\pi_1(D)=\pi_1(\conj D)$ be the group of deck transformations of
the universal covering $\conj\pi:\conj Y\to\conj D$. There exists a representation
$\rho:\Gamma\to\Aut(B)$ such that
$$
\rho(\gamma)\circ\conj F\,(x)=\conj F\circ\gamma\,(x)
\quad\text{for all } x\in\conj Y\text{ and } \gamma\in\Gamma.
$$
\end{lem}

\begin{proof}
The existence of the representation $\rho$ such that the required relation
holds for all $x\in\partial Y$ is a direct consequence of the Poincar\'e--Alexander
theorem~\cite{Al} and was observed by Burns--Shnider~\cite[\S 1]{bs}.
The extension to the entire $\conj Y$ follows by the uniqueness theorem.
\end{proof}

Examples in~\cite{bs} show that the inverse germ 
$_{p'}{\bf f}^{-1}:S\to\partial D$ may not extend along every path in~$S$.
On the other hand, the inverse map does extend along every path in 
the {\it open\/} ball~$B$. To see this, it is enough to prove:

\begin{lem}
\label{ball}
There exists an $\eps>0$ such that every point $x\in Y$ has an open
neighbourhood $V\subset Y$ with the following properties:
\begin{itemize}
\item[(1)] the restriction $F|_V:V\to F(V)$ is biholomorphic,
\item[(2)] $F(V)$ contains the ball of radius $\eps$ centred at $F(x)$
with respect to the Poincar\'e metric on~$B$.
\end{itemize}
\end{lem}

\begin{proof} Let $h$ be the Euclidean metric in $\C^n$. The Poincar\'e 
metric dominates the Euclidean metric in the ball~$B$. In particular, 
the Euclidean ball of radius $R>0$ centred at a point $b\in B$
contains the Poincar\'e ball of the same radius centred at~$b$.

Denote by $\conj F^*h$ the pull-back of the Euclidean metric
to the manifold $\conj Y$. Let further $\Phi\subset Y$ be a fundamental
domain for the action of $\Gamma$ on~$Y$.
Notice that $\Phi$ is a relatively compact subset of~$\conj Y$.
It follows that there exists an $\eps>0$ such that for every
point $x\in\Phi$ the map $\conj F$ is a biholomorphism of the ball of
radius $\eps$ centred at $x$ with respect to the metric $\conj F^*h$
onto the intersection of the Euclidean ball of the same radius
centred at $\conj F(x)$ with the closed ball~$\conj B$. Since the
Euclidean ball about an interiour point of $B$ contains
the Poincar\'e ball with the same centre and radius, we have just
shown that every point in $\Phi$ does indeed have a neighbourhood
with properties (1) and~(2).

Now let $x\in Y$ be an arbitrary point. By the definition of a
fundamental domain, there exists a deck transformation $\gamma\in\Gamma$
such that $\gamma(x)\in\Phi$. Let $W$ be the neighbourhood of $\gamma(x)$
constructed above and set $V=\gamma^{-1}(W)$. By Lemma~\ref{equiv}, we have
$$
F=\rho(\gamma)^{-1}\circ F\circ \gamma.
$$
It follows that $F$ is biholomorphic in $V$ if and only if it is
biholomorphic in~$W$. Furthermore, the image $F(V)=\rho(\gamma)^{-1}(F(W))$
contains the Poincar\'e ball of radius $\eps$ about~$F(x)$ because
$F(W)$ contains the ball of this radius about $F(\gamma(x))$ and the
automorphism $\rho(\gamma)^{-1}$ is an isometry of the Poincar\'e metric.
\end{proof}

Since $B$ is simply connected, it follows that the map $F:Y\to B$ of 
open manifolds is biholomorphic. This proves Theorem~\ref{mainA} and 
the `if' part of Theorem~\ref{c1}.

\begin{rem}[Boundary behaviour I] Once $F$ is known to be biholomorphic,
it is easy to see that $\conj F$ is injective on $\partial Y$. Hence,
$\conj F$ is a biholomorphic map from $\conj Y$ onto $\conj B\setminus A$,
where $A=S\setminus \conj F(\partial Y)$ is a closed subset of the unit sphere. 
In other  words, the closure of the domain $D$ is uniformized by 
$\conj B\setminus A$. However, the subset~$A$ depends on~$D$, 
and its structure remains mysterious.
\end{rem}

Let us now prove the `only if' part of Theorem~\ref{c1}. The argument
is completely independent of the rest of the paper.

\begin{prop} If a Stein strictly pseudoconvex domain $D$ with real analytic
boundary is covered by the open unit ball, then its boundary is spherical.
\end{prop}

\begin{proof} Let $\pi: B\to D$ be the covering. Let $q\in \partial D$ be
an arbitrary point and $U$ a coordinate neighbourhood of $q$ such that
$D\cap U$ is simply connected. Then the germ of the map $\pi^{-1}$ extends
to a biholomorphic map $g$ from $U\cap D$ to an open set in $B$.
A standard argument using the Hopf lemma and the
asymptotics of the Poincar\'e metric on~$B$ (see, e.\,g., \cite{pi0}, \cite{FL}
or~\cite{Su}) shows that $g$ extends to $\partial D\cap U$ as a H\"older
continuous map sending $\partial D\cap U$ to the unit sphere. Note that by
the boundary uniqueness theorem the extension to the boundary is not constant.
Hence,  the extension of $g$ to $\partial D$ is in fact smooth by~\cite{PT}.
Finally, by the reflection principle of Lewy and Pinchuk, the map $g$ extends
biholomorphically to a neighbourhood of $q$, which shows that $\partial D$
is spherical.
\end{proof}

The uniformization theorem imposes strong restrictions on the
topology of spherical domains.

\begin{cor}\label{topol}
Let $D$ be a Stein strictly pseudoconvex domain with spherical boundary. 
Then the higher homotopy groups $\pi_k(D)=0$ for all $k\ge 2$. 
If $D$ is not biholomorphic to the ball, then its fundamental group 
is infinite and contains no non-trivial finite subgroups.
\end{cor}

\begin{proof}
Let us give a purely topological proof of this fact.  The higher homotopy 
groups vanish because the universal covering of $D$ is contractible. 
Suppose that $\pi_1(D)$ contains a non-trivial element of finite prime 
order~$p>0$. Let $\widetilde D$ be the covering of $D$ corresponding
to the subgroup generated by this element. Then $\pi_1(\widetilde D)=\Z/p\Z$ 
and $\pi_k(\widetilde D)=0$ for all $k\ge 2$. By the Hurewicz--Eilenberg--McLane
theorem, the cohomology groups of $\widetilde D$ with any coefficient ring
are isomorphic to the cohomology groups of its fundamental group.
However, $H^k(\Z/p\Z;\Z/p\Z)\ne 0$ for all $k\ge 0$ 
(see, e.\,g., \cite[p.~28]{Ev}) whereas the space $\widetilde D$, which 
is a finite-dimensional manifold, cannot have non-trivial cohomology 
in all positive dimensions. This contradiction shows that every finite 
subgroup of $\pi_1(D)$ is trivial.
\end{proof}

\begin{rem}[Boundary behaviour II]
The Lefschetz theorem for Stein manifolds~\cite{AF},\cite{Bo} tells us
that the homomorphism $\pi_k(\partial D)\to\pi_k(D)$ is an isomorphism
for all $k\le n-2$ and a surjection for $k=n-1$, where $n=\dim_\C D$.
In particular, $\pi_1(\partial D)$ surjects onto $\pi_1(D)$ for all $n\ge 2$,
and we retrieve the result of Burns--Shnider~\cite{bs} that a compact
spherical hypersurface bounding a Stein domain other than the ball
must have infinite fundamental group.

If the complex  dimension~$n\ge 3$, then $\pi_1(\partial D)=\pi_1(D)$.
It follows that, firstly, $\pi_1(\partial D)$ does not possess non-trivial
finite subgroups and, secondly, the covering $\partial Y\to\partial D$
is the universal covering of the boundary. On the other hand, the
paper~\cite{GKL} provides examples of strictly pseudoconvex Stein quotients
of the unit ball in~$\C^2$ having torsion elements in~$\pi_1(\partial D)$.
The same examples show that the covering $\partial Y\to\partial D$ is not
in general the universal covering of the boundary in the two-dimensional case.
\end{rem}

\subsection{Proof of Theorem \ref{mainC}}
The argument follows the familiar pattern. Let $_p{\bf f}:M\to M'$
be the germ of a local equivalence between two Levi non-degenerate
pseudoconcave real analytic compact hypersurfaces in $\CP^n$.
If $M'$ is real algebraic, it follows from Theorem~\ref{Shaf}
that this germ extends as a locally biholomorphic map along
any path in~$M$. Therefore, it extends as a locally biholomorphic
map sending $M$ to $M'$ along any path in a neighbourhood $U\supset M$
provided by Lemma~\ref{LocExtCP}. This extension defines a
locally biholomorphic map $F:\widehat U\to\CP^n$ from the
universal covering of the neighbourhood~$U$.

Now we are in the situation discussed in Example~\ref{Concave}.
In particular, we know that every holomorphic function on $\widehat U$
is constant. Hence, Corollary~\ref{CPRig} shows that the map
$F:\widehat U\to\CP^n$ is the pull-back of an automorphism $L\in\Aut(\CP^n)$.
But this is equivalent to saying that $_p{\bf f}={}_p{\bf L}$,
which completes the proof.\qed

\begin{exm} Let us briefly outline a construction showing that
there exist many topologically different examples of real algebraic 
pseudoconcave hypersurfaces in $\CP^n$. Note that it is enough to exhibit
smooth Levi non-degenerate hypersurfaces in $\CP^n$ because they
can be approximated by real algebraic hypersurfaces, which will
then have the same signature of the Levi form and the same
topology.

Let $Z\subset\CP^n$ be a complex submanifold of dimension~$k\ge 0$.
A well-known result going back to Grauert states that, since the
normal bundle of $Z$ is positive, the boundary $M=\partial U$
of an appropriately chosen tubular neighbourhood $U\supset Z$ has the Levi
form of signature~$(n-k-1,k)$. Here the $n-k-1$ positive directions
are `perpendicular' to $Z$ and the $k$ negative directions are `parallel'
to~$Z$. For instance, a point $(k=0)$ has a strictly pseudoconvex
neighbourhood, and a complex hypersurface $(k=n-1)$ a strictly
pseudoconcave one. The standard real hyperquadrics can be obtained
by this construction from linear subspaces $Z\subset\CP^n$ of appropriate
(co)dimension.

Topologically, $M$ is an $S^{2(n-k)-1}$-bundle over~$Z$.
Thus, the fundamental group of $M$ is isomorphic to that of~$Z$
if the complex codimension of $Z$ is at least~$2$. For instance,
if $Z$ is a complex curve of genus $g>0$ embedded in~$\CP^3$,
then $M$ is a pseudoconcave real hypersurface of Levi signature~$(1,1)$
with infinite fundamental group.
\end{exm}


\end{document}